\documentclass[12pt,reqno]{amsart}
\usepackage{graphicx, amssymb, color,pdfsync}
\usepackage[all,cmtip]{xy}
\numberwithin{equation}{section}
\usepackage{mathrsfs}
\usepackage{hyperref}

\hypersetup{
    colorlinks=true,       
    linkcolor=blue,          
    citecolor=blue,        
    filecolor=blue,      
    urlcolor=blue           
}
\usepackage{tikz}
\usetikzlibrary{matrix,arrows}
\DeclareMathOperator{\Aut}{Aut}
\usepackage[switch]{lineno}
\usepackage{yfonts}
\advance\hoffset by -0.9 truecm

\begin{document}
\newcommand{\s}{\vspace{0.2cm}}

\newtheorem{theo}{Theorem}
\newtheorem{prop}{Proposition}
\newtheorem{coro}{Corollary}
\newtheorem{lemm}{Lemma}
\newtheorem{remas}{Remark}

\newtheorem{claim}{Claim}
\newtheorem{example}{Example}
\theoremstyle{remark}
\newtheorem*{rema}{\it Remarks}
\newtheorem*{rema1}{\it Remark}
\newtheorem*{defi}{\it Definition}
\newtheorem*{theo*}{\bf Theorem}
\newtheorem*{coro*}{Corollary}

\title[On $pq$-fold regular covers of the projective line]{On $pq$-fold regular covers\\ of the projective line}
\date{}

\author{Sebasti\'an Reyes-Carocca}
\address{Departamento de Matem\'atica y Estad\'istica, Universidad de La Frontera, Avenida Francisco Salazar 01145, Temuco, Chile.}
\email{sebastian.reyes@ufrontera.cl}

\thanks{Partially supported by Fondecyt Grants 11180024, 1190991 and Redes Grant 2017-170071}

\begin{abstract} Let $p$ and $q$ be odd prime numbers.
In this paper we study non-abelian $pq$-fold regular covers of the projective line, determine algebraic models for some special cases  and provide a general isogeny decomposition of the corresponding Jacobian varieties. We also give a classification and description of the one-dimensional families of compact Riemann surfaces as before. \end{abstract}
\maketitle
\thispagestyle{empty}

\section{Introduction and statement of the results}

Compact Riemann surfaces (or, equivalently, smooth  complex  projective algebraic curves) and their automorphism groups have been extensively studied since the nineteenth century.  Foundational results concerning that are: \begin{enumerate}
\item if the genus of the compact Riemann surface is greater than one
 then its automorphism group is finite (see \cite{Hu} and \cite{SS}, and also \cite{FK}), and 
\item each finite group acts as a group of automorphisms of some compact Riemann surface of a suitable genus greater than one (see \cite{Gre} and also \cite{Ga}). \end{enumerate}

\s

A general problem that arises naturally with regard to this is to determine necessary and sufficient conditions under which a  given group  acts as a group of automorphisms of a compact Riemann surface satisfying some prescribed conditions. This problem was successfully studied for cyclic groups by Harvey in \cite{H2} and soon after for abelian groups by Maclachlan in \cite{Mcl}. The same problem for dihedral groups, among other aspects, was completely solved by Bujalance, Cirre, Gamboa and Gromadzki in \cite{buja}. See also \cite{die}.

\s

This article is devoted to study those compact Riemann surfaces  that are branched $pq$-fold regular covers of the projective line, where $p$ and $q$ are prime numbers. Since the abelian and dihedral cases have been already classified, we shall consider compact Riemann surfaces endowed with a non-abelian group of automorphisms  isomorphic to the semidirect product of two cyclic groups of odd prime order, in such a way that the corresponding orbit space is isomorphic to  the projective line. 

\s

Let $p$ and $q$ be odd primes such that $p$ divides $q-1$ and let $r$ be a  primitive $p$-th root of unity  in the field of $q$ elements.  Throughout the article the unique non-abelian group of order $pq$ will be denote by $$G_{p,q}:=\langle a,b : a^q=b^p=1, bab^{-1}=a^r\rangle \cong C_q \rtimes C_p.$$

The first result of this paper establishes a simple necessary and sufficient condition for $G_{p,q}$ to act on a compact Riemann surface of genus greater than one. In order to state it, we need to bring in the concept of {\it signature}. The tuple $$(\gamma; k_1, \ldots, k_l) \in \mathbb{Z}^{l+1}  \mbox{ where }  \gamma \geqslant 0 \mbox{ and } k_i \geqslant 2$$is called the {\it signature} of the action of a group $G$ on a compact Riemann surface $S$ if the genus of the orbit space $S/G$ is $\gamma$ and the  branched regular covering map $$S \to S/G$$ ramifies over exactly $l$ values $y_1, \ldots, y_l$ and the fiber over $y_i$ consists of points with $G$-isotropy group of order $k_i$ for each $i \in \{1, \ldots, l\}.$ Note that $k_i$ divides $|G|.$

\begin{theo} \label{tepi}
Let $n,m \geqslant 0$ be integers such that $n+m \geqslant 3$. There exists a compact Riemann surface $S$ of genus greater than one endowed with a group of automorphisms isomorphic to $G_{p,q}$ acting on it with signature $$s_{n,m}:=(0; p, \stackrel{n}{\ldots}, p, q, \stackrel{m}{\ldots}, q)$$if and only if $n \geqslant 2.$ In this case, the genus of $S$ is $g=1-pq+nq(\tfrac{p-1}{2})+mp(\tfrac{q-1}{2}).$ 
\end{theo}

Let $\mathscr{M}_g$ denote the moduli space of compact Riemann surfaces of genus $g \geqslant 2.$ For each pair of integers $n, m$ as in Theorem \ref{tepi}, the set of compact Riemann surfaces admitting a group of automorphisms isomorphic to $G_{p,q}$ acting on them with signature $s_{n,m}$ form a {\it family} of complex dimension $n+m-3$ in the singular locus of $\mathscr{M}_g$  (see Subsection \S\ref{sta} for a precise definition of {\it family}). 

\s

{\bf Notation.} We shall denote the above introduced family by $\mathscr{C}_{n,m}.$

\subsection*{ \it The zero-dimensional case} A compact Riemann surface is called {\it quasiplatonic} if it has a group of automorphisms such the signature of the action is of the form $$(0; k_1, k_2, k_3).$$

Observe that among the compact Riemann surfaces of Theorem \ref{tepi}, the quasiplatonic ones split into two cases: those with signature $$s_{2,1}=(0; p,p,q) \,\, \mbox{ and }\,\, s_{3,0}=(0; p,p,p).$$ 

These compact Riemann surfaces form two zero-dimensional families; namely, there are, up to isomorphism, finitely many of them. Streit and Wolfart in \cite{SW} succeeded in describing these Riemann surfaces in full detail. More precisely, they determined algebraic models, the number of isomorphism classes, their automorphism groups and their minimal fields of definition, among other aspects. 

\s

Given a compact Riemann surface $S$ with a group of automorphisms $G$, a natural question that arises is to decide whether or not $S$ admits more automorphisms. This  is a challenging problem and its answer depends not only on the signature $s$ of the action of $G$, but also on the geometry of a  fundamental domain for the  surface $S$. 

Singerman in \cite{singerman2} determined all those possible signatures $s$ for which the Riemann surface might have more automorphisms (see Subsection \S \ref{didid}). A  direct consequence of his results is that if the complex dimension of the family $\mathscr{C}_{n,m}$ is greater than one, then $G_{p,q}$ is the full automorphism group of the surfaces lying in the interior of $\mathscr{C}_{n,m}$. The one-dimensional case will be considered later in this paper. 
\s

We recall the obvious fact that $G_{p,q}$ has a unique normal subgroup $N$ of order $q,$ and $q$ pairwise conjugate subgroups of order $p;$ let $H$ be one of them. All along the paper $N$ and $H$ will be used to denote these subgroups. For each $S \in \mathscr{C}_{n,m}$ we consider the   regular covering maps $$S \to X:=S/N \,\, \mbox{ and } \,\, S \to Y:= S/H $$given, respectively, by the action  of $N$ and $H$ on $S$.  As an application of the Riemann-Hurwitz formula, we see that the genera of $X$ and $Y$ are
$$g_X=\tfrac{p-1}{2}(n-2) \, \,\mbox{ and } \, \,\, g_Y=\tfrac{q-1}{2}(m-2)+\tfrac{p-1}{2}\tfrac{q-1}{p}n.$$

Although the literature still shows few general results in this direction, there is a great interest in providing explicit  descriptions of compact Riemann surfaces as algebraic curves. See, for instance, \cite{M1}, \cite{M2}, \cite{M3}, \cite{C1}, \cite{Shaska}, \cite{yojpaa}, \cite{yoanita2},  \cite{T1} and \cite{T2}.

The following result provides an algebraic description for each compact Riemann surface $S$ lying in the family $\mathscr{C}_{2,m}$ for  $m \geqslant 1.$ This description is given in terms of a singular plane algebraic curve  whose projective desingularization is isomorphic to $S.$

\begin{theo} \label{modelo}
Let $n,m \geqslant 0$ be integers such that $n+m \geqslant 3.$
\begin{enumerate} 
\s
\item Let $S$ belong to the interior of the family $\mathscr{C}_{n,m}.$ Then $S$ is cyclic $\hat{q}$-gonal for some odd prime $\hat{q}$ if and only if $n=2$ and $\hat{q}=q.$
\s

\item An affine singular algebraic model for $S \in \mathscr{C}_{2,m}$ in $\mathbb{C}^2$ is $$y^q=\Pi_{i=0}^{p-1}(x-\omega^i)^{r^i}  \Pi_{k=2}^{m}\Pi_{i=0}^{p-1}(x- \lambda_k \omega^i)^{r^i}$$where $\omega$ is a primitive $p$-th root of unity and $\lambda_2, \ldots, \lambda_{m}$ are non-zero complex numbers such that $\lambda_k^p \neq 1$ for each $k$ and $\lambda_k \notin \{ \lambda_j \omega^i : i \in \mathbb{Z}_p\}$ for all $j \neq k.$ 

\s

\item In the previous model, the group of automorphisms of $S$ isomorphic to $G_{p,q}$ is generated by $$A(x,y)=(x, \xi y) \,\, \mbox{ and }\,\, B(x,y)=(\omega x, \varphi(x)y^r)$$where $\xi$ is a primitive $q$-th root of unity,  $$\varphi(x)=\omega^{me} [  (x-\omega^{p-1})\Pi_{k=2}^m(x-\lambda_k\omega^{p-1})]^{e(1-r)}$$and $e \in \mathbb{Z}$ is chosen to satisfy $1+r+\cdots + r^{p-1}=e q.$
\end{enumerate}
\end{theo}

\subsection*{ \it The one-dimensional case} Observe that, by Theorem \ref{tepi}, there exist three complex one-dimensional families of compact Riemann surfaces that are non-abelian $pq$-fold branched regular covers of the projective line. More precisely, these families are:

\s
\s

\begin{center}
\begin{tabular}{|c|c|c|}  
\hline

\, family \, &  genus  &  signature   \\ [0.4ex]\hline 
$\mathscr{C}_{2,2}$  & $(p-1)(q-1)$ & $\, s_{2,2}=(0; p,p,q,q) \,$ \\ [0.4ex]
$\mathscr{C}_{4,0}$ &  $1+q(p-2)$ & $s_{4,0}=(0; p,p,p,p)$  \\[0.4ex] 
$\mathscr{C}_{3,1}$  & $\, 1+pq-\tfrac{1}{2}(p+3q) \,$ & $s_{3,1}=(0; p,p,p,q)$ \\[0.4ex] 
\hline
\end{tabular}
\end{center}

\s
\s
The following three theorems provide a description of these  families. More precisely, we give upper bounds for the number of equisymmetric strata they consists of (see Subsection \S\ref{sta} for a precise definition of {\it equisymmetric stratum}), determine their automorphism groups and show that the complement of the interior in the closure of each of them is non-empty.

\begin{theo} \label{familia22}
\mbox{}
\begin{enumerate}

\item The family $\mathscr{C}_{2,2}$ consists of at most $q(p-1)$ equisymmetric strata.
\s
\item If $S$ belongs to the interior of $\mathscr{C}_{2,2}$ then the automorphism group of $S$ is isomorphic to $G_{p,q}.$ 
\s
\item The closure of the family $\mathscr{C}_{2,2}$ contains a quasiplatonic Riemann surface with a group of automorphisms of order $2pq$ isomorphic to $$C_q \rtimes C_{2p} \,\mbox { acting with signature } \,(0; q, 2p, 2p).$$
 
\end{enumerate}
\end{theo}

\begin{theo} \label{familia40} \mbox{}
\begin{enumerate}

\item  The family $\mathscr{C}_{4,0}$ consists of at most $q(p-1)(p^2-3p+3)$ equisymmetric strata.

\s

\item If $S$ belongs to the interior of $\mathscr{C}_{4,0}$ then either  the automorphism group of $S$ is isomorphic to $G_{p,q}$ or  to $C_q \rtimes C_{2p}$ acting on it with signature $(0; 2,2,p,p).$
\s
\item The closure of the family $\mathscr{C}_{4,0}$ contains  a quasiplatonic Riemann surface with a group of automorphisms of order $4pq$ isomorphic to $$(C_q \rtimes C_{2p})\rtimes C_2 \,\mbox { acting with signature } \,(0; 2, 2p, 2p).$$

\item The compact Riemann surface $X=S/N$ is hyperelliptic of genus $p-1,$ and an affine singular algebraic model for it in $\mathbb{C}^2$ is $$y^2 =(x^p-1)(x^p-\mu^p)$$where $\mu$ is a non-zero complex number such that $\mu^p \neq 1.$ 
\end{enumerate}
\end{theo}

\begin{theo} \label{familia31}
\mbox{}
\begin{enumerate}

\item The family $\mathscr{C}_{3,1}$ consists of at most $q(p-1)(p-2)$ equisymmetric strata.
\s

\item If $S$ belongs to the interior of $\mathscr{C}_{3,1}$ then the automorphism group of  $S$ is isomorphic to $G_{p,q}.$ 
\s

\item The closure of the family $\mathscr{C}_{3,1}$ contains a  quasiplatonic Riemann surface with a group of automorphisms of order $2pq$ isomorphic to $$(C_q \rtimes C_{p}) \times C_2 \,\mbox { acting with signature } \,(0; p, 2q, 2q).$$
\end{enumerate}
\end{theo}

\s

\begin{rema}
\begin{enumerate} \mbox{}
\item If $S$ belongs to the interior of $ \mathscr{C}_{n,m}$ then $S$ is non-hyperelliptic.
\item The family $\mathscr{C}_{4,0}$ with $p=3$ was studied in \cite{IRC}; see also \cite{IJRC}.
\item The algebraic description of $S \in \mathscr{C}_{2,m}$ extends the ones in \cite{SW}; see also \cite{CM}.

\end{enumerate}
\end{rema}

\subsection*{ \it Jacobian varieties} Let $S$ be a compact Riemann surface of genus $g.$ We denote by $$\mathscr{H}^1(S, \mathbb{C})^* \,\, \mbox{ and }\,\, H_1(S, \mathbb{Z})$$ the dual of the complex vector space of dimension $g$ of its 1-forms and  its first integral homology group, respectively. The {\it Jacobian variety} of $S,$ defined as the quotient $$JS=\mathscr{H}^1(S, \mathbb{C})^*/H_1(S, \mathbb{Z}),$$is an irreducible principally polarized abelian variety of dimension $g.$ The importance of  the Jacobian variety of a compact Riemann surface lies, partially, in the classical Torelli's  theorem, which ensures that, up to isomorphism, the Riemann surface is uniquely determined by its Jacobian variety. Namely, $$S_1 \cong S_2 \,\, \,\mbox{ if and only if }\,\,\, JS_1 \cong JS_2.$$

Is is known that the action of a group $G$ on $S$  induces a isogeny decomposition  $$JS \sim B_1^{e_1} \times \cdots \times  B_s^{e_s} $$in terms of abelian subvarieties $B_1, \ldots, B_s$ of $JS$ in such a way that $$A_i=B_i^{e_i}  \mbox{ for } i=1, \ldots, s$$are pairwise non-$G$-isogenous.  

\s

It is worth emphasizing that whereas this decomposition only depends on the algebraic structure of the group, the dimension of the factors $B_i$ does depend on the way the group acts. More precisely, it depends on the signature  and, in addition, on the equisymmetric stratum to which $S$ belongs. An example of a family exhibiting in an explicit manner this dependence can be found in \cite[Remark 2(6)]{familias}.

\begin{theo} \label{jaco} 
Let $S$ be a compact Riemann surface lying in the family $\mathscr{C}_{n,m}.$ Then the Jacobian variety $JS$ of $S$ decomposes, up to isogeny, as the product $$JS \sim B_1 \times B_2^p$$where $B_1$ and $B_2$ are abelian subvarieties of $JS$ of dimension$$\tfrac{p-1}{2}(n-2) \, \,\mbox{ and } \, \,\, \tfrac{q-1}{2}(m-2)+\tfrac{p-1}{2}\tfrac{q-1}{p}n.$$
\end{theo}

We point out that in the previous theorem, the dimension of the factors $B_1$ and $B_2$ only depends on the signature of the action.

\s

Let $K$ be a group of automorphisms of a compact Riemann surface $S$ and let \begin{equation} \label{holass}S \to R=S/K\end{equation} be the associated regular covering map. Assume the genus of $R$ to be $\gamma \geqslant 1.$  Then \eqref{holass} induces a homomorphism between the associated Jacobian varieties\begin{equation*}\label{indi}JR \to JS\end{equation*}whose image  is an abelian subvariety of $JS$ of dimension $\gamma$ which is isogenous to $JR$. 

\s

Keeping the same notations as before, we observe that$$g_X=\dim(JX)=\dim(B_1) \,\, \mbox{ and } \,\, g_Y=\dim(JY)=\dim(B_2).$$

As a matter of fact,  the following result holds.
\begin{theo} \label{prop3}
Let $S$ be a compact Riemann surface lying in the family $\mathscr{C}_{n,m}.$ With the notations of Theorem \ref{jaco}, there exist isogenies $$JX \to B_1 \,\, \mbox{ and }\,\, JY \to B_2$$and, in particular, $JS$ is isogenous to the product of Jacobians of quotients of $S$  $$JS \sim JX \times JY^p$$
\end{theo}
\s

This article is organized as follows. Section \S \ref{preli} will be devoted to briefly review the basic background: Fuchsian groups and group actions on compact Riemann surfaces and Jacobian varieties. The proof of the theorems are given in the remaining sections. 

\section{Preliminaries} \label{preli}

\subsection{Fuchsian groups}  A {\it Fuchsian group} is a discrete group of automorphisms of $$\mathbb{H}=\{z \in \mathbb{C}: \mbox{Im}(z) >0 \}.$$  

If $\Delta$ is a Fuchsian group and the orbit space $\mathbb{H}/{\Delta}$ given by the action of $\Delta$ on $\mathbb{H}$ is  compact, then the algebraic structure of $\Delta$ is determined by its {\it signature}: \begin{equation} \label{sig} s(\Delta)=(\gamma; k_1, \ldots, k_l),\end{equation}where  the genus of   $\mathbb{H}/{\Delta}$ is $\gamma$ and $k_1, \ldots, k_l$ are the branch indices in the universal canonical projection $\mathbb{H} \to \mathbb{H}/{\Delta}.$ 

\s

If  $\Delta$ is a Fuchsian group of signature \eqref{sig} then $\Delta$ has a canonical presentation in terms of generators $\alpha_1, \ldots, \alpha_{\gamma}$, $\beta_1, \ldots, \beta_{\gamma},$ $ x_1, \ldots , x_l$ and relations
\begin{equation}\label{prese}x_1^{k_1}=\cdots =x_l^{k_l}=\Pi_{i=1}^{\gamma}[\alpha_i, \beta_i] \Pi_{i=1}^l x_i=1,\end{equation}where the brackets stands for the commutator, and the hyperbolic area of each fundamental region of $\Delta$ is given by $$\mu(\Delta)=2 \pi [2\gamma-2 + \Sigma_{i=1}^l(1-\tfrac{1}{k_i})].$$

Let $\Delta'$ be a group of automorphisms of $\mathbb{H}.$ If $\Delta$ is a subgroup of $\Delta'$ of finite index, then $\Delta'$ is a Fuchsian group and their hyperbolic areas are related by the Riemann-Hurwitz formula $$\mu(\Delta)= [\Delta' : \Delta] \cdot \mu(\Delta').$$
 
The Teichm\"{u}ller space of $\Delta$ is a complex analytic manifold homeomorphic to the complex ball of dimension $3\gamma-3+l$. See, for instance, \cite{singerman2} for further details.
\subsection{Group actions on Riemann surfaces} \label{didid}Let $S$ be a compact Riemann surface of genus $g \geqslant 2$ and let $\mbox{Aut}(S)$ denote its  automorphism group. A finite group $G$ acts on $S$ if there is a group monomorphism $G\to \Aut(S).$ The orbit space $S/G$ of the action of $G$ on $S$ is naturally endowed with a compact Riemann surface structure such that the canonical projection $S \to S/G$ is holomorphic. 

\s

By the classical uniformization theorem, there is a unique, up to conjugation, Fuchsian group $\Gamma$ of signature $(g; -)$ such that $S \cong \mathbb{H}/{\Gamma}.$ Moreover, $G$ acts on $S$ if and only if there is a Fuchsian group $\Delta$  together with a group  epimorphism \begin{equation*}\label{epi}\theta: \Delta \to G \, \, \mbox{ such that }  \, \, \mbox{ker}(\theta)=\Gamma.\end{equation*}

 The action is said to be represented by {\it the surface-kernel epimorphism} $\theta$; henceforth, we write {\it ske} for short. It is said that    $G$ acts on $S$ with signature $s(\Delta).$ Note that this definition agrees with the one given in the introduction.
 
\s

Assume that $G'$ is a finite group such that $G \leqslant G'.$ Then the action of $G$ on $S$ represented by the ske $\theta$ is said to {\it extend} to an action of $G'$ on $S$ if:\begin{enumerate}
\item there is a Fuchsian group $\Delta'$ containing $\Delta,$ 
\item the Teichm\"{u}ller spaces of $\Delta$ and $\Delta'$ have the same dimension, and
\item there exists a ske  $$\Theta: \Delta' \to G' \, \, \mbox{ in such a way that }  \, \, \Theta|_{\Delta}=\theta   \mbox{ and } \mbox{ker}(\theta)=\mbox{ker}(\Theta).$$
\end{enumerate} 

{\it Maximal actions} are those that  cannot be extended in the previous sense. A complete list of pairs of signatures of Fuchsian groups $\Delta$ and $\Delta'$ for which it may be possible to have an extension as before was provided by Singerman in \cite{singerman2}. 

\subsection{Equivalence of actions} \label{strati} Let $\text{Hom}^+(S)$ denote the group of orientation preserving self-homeomorphisms of $S.$ Two actions $\psi_i: G \to \mbox{Aut}(S)$ of $G$ on $S$ with $i=1,2$ are  {\it topologically equivalent} if there exist $\omega \in \Aut(G)$ and $f \in \text{Hom}^+(S)$ such that
\begin{equation}\label{equivalentactions}
\psi_2(g) = f \psi_1(\omega(g)) f^{-1} \hspace{0.5 cm} \mbox{for all} \,\, g\in G.
\end{equation}

Each $f$ satisfying \eqref{equivalentactions} yields an automorphism $f^*$ of $\Delta$ where $\mathbb{H}/{\Delta} \cong S/G$. If $\mathscr{B}$ is the subgroup of $\mbox{Aut}(\Delta)$ consisting of them, then $\mbox{Aut}(G) \times \mathscr{B}$ acts on the set of skes defining actions of $G$ on $S$ with signature $s(\Delta)$ by $$((\omega, f^*), \theta) \mapsto \omega \circ \theta \circ (f^*)^{-1}.$$  

Two skes $\theta_1, \theta_2 : \Delta \to G$ define topologically equivalent actions if and only if they belong to the same $(\mbox{Aut}(G) \times \mathscr{B})$-orbit (see \cite{Brou},  \cite{H} and \cite{McB}). If the genus of $S/G$ is zero then $\mathscr{B}$ is generated by the {\it braid transformations}  $\Phi_{i}$, for $1 \leqslant i  < l,$ defined by \begin{equation*} \label{braid} x_i \mapsto x_{i+1}, \hspace{0.3 cm}x_{i+1} \mapsto x_{i+1}^{-1}x_{i}x_{i+1} \hspace{0.3 cm} \mbox{ and }\hspace{0.3 cm} x_j \mapsto x_j \mbox{ when }j \neq i, i+1.\end{equation*}

\subsection{Stratification of the moduli space and families} \label{sta} We denote the Teichm\"{u}ller space of a Fuchsian group of signature $(g; -)$ by $T_g$. It is well-known that the moduli space $\mathscr{M}_g$ arises as the quotient space$$\pi : T_g \to \mathscr{M}_g:=T_g/\mbox{Mod}_g$$given by the action of the mapping class group $\mbox{Mod}_g$ of genus $g$ on $T_g.$ Observe that $\mathscr{M}_g$ is endowed with the quotient topology induced by $\pi$. Moreover, $\mathscr{M}_g$ is endowed with a structure of complex analytic space of dimension $3g-3,$ and for $g \geqslant4$ its singular locus agrees with the branch locus of $\pi$ and correspond to set of points representing compact Riemann surfaces with non-trivial automorphisms. See, for instance, \cite[Section 2]{b}  for further details.

\s

Following \cite[Theorem 2.1]{b}, the singular locus of $\mathscr{M}_g$ admits an {\it equisymmetric stratification} $$\mbox{Sing}(\mathscr{M}_g)= \cup_{G, \theta} \bar{\mathscr{M}}_g^{G, \theta}$$ where 
each {\it equisymmetric stratum} $\mathscr{M}_g^{G, \theta}$, if nonempty, corresponds to one topological class of maximal actions (see also \cite{H}). More precisely:

\begin{enumerate}

\item  ${\mathscr{M}}_g^{G, \theta}$ consists of surfaces of genus $g$ with automorphism group isomorphic to $G$ such that the action is topologically equivalent to $\theta$,

\item the {\it closure} $\bar{\mathscr{M}}_g^{G, \theta}$ of  ${\mathscr{M}}_g^{G, \theta}$ is a closed irreducible algebraic subvariety of $\mathscr{M}_g$ and consists of surfaces  of genus $g$ with a group of automorphisms isomorphic to $G$ such that the action is  topologically equivalent to $\theta$, and

\item  if ${\mathscr{M}}_g^{G, \theta}$ is nonempty then it is a smooth, connected,
locally closed algebraic subvariety of $\mathscr{M}_{g}$ which is Zariski dense in
$\bar{\mathscr{M}}_g^{G, \theta}.$ 
\end{enumerate}

\s

The aforementioned stratification was the key ingredient in finding all those values of $g$ for which the singular locus of $\mathscr{M}_g$ is connected; see, for instance, \cite{BCI} and \cite{BCIP}.

\s

As mentioned in the introduction, we shall employ the following terminology.

\begin{defi} Let $G$ be a group and let $s$ be a signature. The subset of $\mathscr{M}_g$ of  all those  compact Riemann $S$ surfaces of genus $g$  with a group of automorphisms isomorphic to $G$ acting with signature $s$ will be called a {\it closed family} or simply a {\it family}.
\end{defi}

Note that  the {\it interior} of the family consists of those Riemann surfaces whose  automorphism group is isomorphic to $G$ and  is formed by finitely many equisymmetric strata which are in correspondence with the pairwise non-equivalent topological actions of $G.$  In addition, the {\it closure} of the family is formed by those surfaces whose automorphism group contains $G$. If the signature of the action of $G$ on $S$ is \eqref{sig} then the dimension of the family is $3\gamma-3+l.$ 

\subsection{Decomposition of Jacobians} \label{ii3} It is  well-known that if $G$ acts on a compact Riemann surface $S$ then this action  induces a $\mathbb{Q}$-algebra homomorphism $$\Phi : \mathbb{Q} [G] \to \mbox{End}_{\mathbb{Q}}(JS)=\mbox{End}(JS) \otimes_{\mathbb{Z}} \mathbb{Q},$$from the rational group algebra of $G$ to the rational endomorphism algebra of $JS.$

For each $ \alpha \in {\mathbb Q}[G]$ we define the abelian subvariety $$A_{\alpha} := {\textup Im} (\alpha)=\Phi (n\alpha)(S) \subset JS$$where $n$ is some positive integer chosen in such a way that $n\alpha \in {\mathbb Z}[G]$.

 Let  $W_1, \ldots, W_s$ be the rational irreducible representations of $G$. For each $W_j$ we denote by $V_j$ a complex irreducible representation of $G$ associated to it.  As proved in \cite[Theorem 2.2]{l-r} (see also \cite[Section 5]{cr}) the decomposition $1=e_1 + \cdots + e_s,$ where $e_j  \in \mathbb{Q}[G]$ is a central idempotent computed explicitly from $W_j$, yields an isogeny $$JS \sim A_{e_1} \times \cdots \times A_{e_s}$$
which is $G$-equivariant.  Moreover, there are idempotents $f_{j1},\dots, f_{jn_j}$ such that $$e_j=f_{j1}+\dots +f_{jn_j} \,\, \mbox{ where } \,\, n_j=d_j/s_j$$is the quotient of the degree $d_j$  and the Schur index $s_j$ of $V_j$. If $B_j=A_{f_j1}$ then we have the isogeny decomposition 
\begin{equation*} \label{m1}
JS \sim B_1^{n_1} \times \cdots \times B_s^{n_s} 
\end{equation*}
called the {\it group algebra decomposition} of $JS$ with respect to $G$. See  also \cite{RCR}.

\s

Assume that \eqref{sig} is the signature of the action of $G$ on $S$  and that this action is represented by $\theta: \Delta \to G,$ with $\Delta$ as in \eqref{prese}. If $W_1$ is the trivial representation, then  $\dim(B_1)=\gamma$ and $n_1=1.$ If $2 \leqslant j \leqslant s$ then, 
 following \cite[Theorem 5.12]{yoibero} 
\begin{equation}\label{ooii}
\dim (B_{j})=c_j[d_j(\gamma -1)+\tfrac{1}{2}\Sigma_{i=1}^l (d_j-d_{j}^{\langle \theta(x_i) \rangle} )]\end{equation} where $c_j$ is the degree of the extension $\mathbb{Q} \le L_j$ with $L_{j}$ denoting a minimal field of definition for $V_j.$
\s

The problem of decomposing Jacobian varieties is old and goes back to Wirtinger \cite{Wir} and Schottky and Jung \cite{SJ}. For recent works concerning that for special classes of groups we refer to the articles     \cite{d1}, \cite{Do}, \cite{FP}, \cite{nos}, \cite{IJR}, \cite{PA}, \cite{ReRu}, \cite{israel}, \cite{Ri} and \cite{Sa}.

\section{Proof of Theorem \ref{tepi}} \label{s3}
 Let $n,m \geqslant 0$ be integers such that $n+m \geqslant 3$ and let $\Delta$ be a Fuchsian group of signature $s_{n,m}=(0; p, \stackrel{n}{\ldots}, p, q, \stackrel{m}{\ldots}, q)$ with canonical presentation $$\langle x_1, \ldots, x_n, y_1, \ldots, y_m : x_1^p=\cdots=x_n^p=y_1^q=\cdots=y_m^q=\Pi_{i=1}^nx_i \Pi_{j=1}^my_j=1\rangle.$$ 
 
Assume the existence of a compact Riemann surface $S$  with a group of automorphisms isomorphic to $G_{p,q}$ acting on it with signature $s_{n,m}.$ Then there exists a ske $$\theta: \Delta \to G_{p,q} \,\, \mbox{ such that }\,\, S \cong \mathbb{H}/\mbox{ker}(\theta).$$
 
Observe that  there is no homomorphism  $\Delta \to G_{p,q}$ provided that $n=1,$ since otherwise the image of $x_1 y_1 \cdots y_m$ would not be trivial. In addition, every homomorphism $\Delta \to G_{p,q}$ with $n=0$ is non-surjective, since $b$ does not belong to its image. Hence, if there exists a compact Riemann surface $S$ as before  then $n \geqslant 2.$ 
 
 \s
 
 Conversely, for each integer $n \geqslant 2$  we shall construct a ske $\theta: \Delta \to G_{p,q}$ explicitly. For the sake of simplicity, we shall identify $\theta$ with the tuple $$\theta=(\theta(x_1), \ldots, \theta(x_n),\theta(y_1), \ldots, \theta(y_m))$$and write $(x,y, \stackrel{t}{\ldots}, x,y)$ to denote the $2t$-uple with $x$ in the $i$-th entry for $i$ odd, and $y$ in the $i$-th entry for $i$ even. 

Assume $m=0.$ If $n$ is even then consider $$(b,b^{-1}, \stackrel{t}{\ldots}, b,b^{-1}, ab, (ab)^{-1})$$where $t=\tfrac{n}{2}-1.$ If $n$ is odd and $n \geqslant 5$ then consider $$(b,b^{-1}, \stackrel{t}{\ldots}, b,b^{-1}, ab, (ab)^{-1}, b^2, b^{-1}, b^{-1})$$where $t=\tfrac{n-5}{2}.$ If $n=3$ consider $(ab,b,(ab^2)^{-1}).$ 

\s

Assume $m=1.$ If $n$ is even then consider $$(b,b^{-1}, \stackrel{t}{\ldots}, b,b^{-1}, b, (ab)^{-1},a)$$where $t=\tfrac{n}{2}-1.$ If $n$ is odd then consider $$(b,b^{-1}, \stackrel{t}{\ldots}, b,b^{-1}, b^2,b^{-1}, (ab)^{-1},a)$$where $t=\tfrac{n-3}{2}.$

\s

Assume $m \geqslant 2.$ If $n$ and $m$ are even then consider $$(b,b^{-1}, \stackrel{t_1}{\ldots}, b,b^{-1}, a,a^{-1}, \stackrel{t_2}{\ldots}, a,a^{-1})$$where $t_1=\tfrac{n}{2}$ and $t_2=\tfrac{m}{2}.$ If $n$ and $m$ are odd then consider $$(b,b^{-1}, \stackrel{t_1}{\ldots}, b,b^{-1}, (b^2, b^{-1}, b^{-1}),  a,a^{-1}, \stackrel{t_2}{\ldots}, a,a^{-1}, (a^2, a^{-1}, a^{-1}))$$where $t_1=\tfrac{n-3}{2}$ and $t_2=\tfrac{m-3}{2}.$ If $n$ is even and $m$ is odd then consider $$(b,b^{-1}, \stackrel{t_1}{\ldots}, b,b^{-1},  a,a^{-1}, \stackrel{t_2}{\ldots}, a,a^{-1}, (a^2, a^{-1}, a^{-1}))$$where $t_1=\tfrac{n}{2}$ and $t_2=\tfrac{m-3}{2}.$ Finally, if $n$ is odd and $m$ is even then consider $$(b,b^{-1}, \stackrel{t_1}{\ldots}, b,b^{-1}, (b^2, b^{-1}, b^{-1}),  a,a^{-1}, \stackrel{t_2}{\ldots}, a,a^{-1})$$where $t_1=\tfrac{n-3}{2}$ and $t_2=\tfrac{m}{2}.$

\s

The value of $g$ is computed as an application of the Riemann-Hurwitz formula to the branched regular covering map $$\mathbb{H}/\mbox{ker}(\theta) \to \mathbb{H}/\Delta$$induced by the inclusion $\mbox{ker}(\theta) \lhd \Delta.$

\section{Proof of Theorem \ref{modelo}}

\subsection*{Proof of statement (1)} The sufficient condition is obvious. Assume that $S$ admits a cyclic $\hat{q}$-gonal morphism for some odd prime number $\hat{q}.$ If  the automorphism group of $S$ is isomorphic to $G_{p,q}$ then $\hat{p}$ equals $p$ or $q.$ The former case is not possible since the genus of $Y$ is always positive. Thus, $\hat{q}=q$ and necessarily $n=2.$ If the automorphism group of $S$ has order strictly greater than $pq$ then, by \cite[Theorem 1]{singerman2},  the pair $(n,m)$ might only belong to$$ \{(2,1), (3,0), (2,2), (4,0)\}.$$However, for the first two cases the automorphism group of $S$ is isomorphic to either $G_{p,q}$ or $C_q \rtimes C_{2p}$ as proved in Theorems 1 and 2 in  \cite{SW}, and the same fact holds for the third and fourth ones as it will be proved later in our Theorems \ref{familia22} and \ref{familia40}. 

\subsection*{Proof of statement (2)} Let $S$ be a compact Riemann surface lying in the family $\mathscr{C}_{2,m}$ with $m \geqslant 1.$  The $q$-gonal morphism $S \to X \cong \mathbb{P}^1$ given by the action of $N=\langle a \rangle$ on $S$ ramifies over $pm$ values, all of them marked with $q.$ If we denote these  values by $$u_{i,k} \,\, \mbox{ for }\,\,i \in \{0, \ldots, p-1\} \mbox{ and } k \in \{1, \ldots, m\}$$and assume that none of them equal $\infty$ then, following \cite{Gab} (see also \cite{H} and \cite{W}), the affine singular algebraic curve $$y^{q}=\Pi_{i=0}^{p-1}\Pi_{k=1}^{m}(x-u_{i,k})^{n_{i, k}} $$is (after normalization) isomorphic to $S,$ for suitable values $1 \leqslant n_{i, k} \leqslant q-1$ in such a way that their sum is congruent to 0 modulo $q.$ Note that $$P= G_{p,q}/N \cong C_p = \langle \beta : \beta^p=1 \rangle$$acts on $X$ with signature $(0; p,p).$ In other words, the action of $P$ on $X$ has two fixed points and $m$ orbits of length $p.$ It follows that, up to a M\"{o}bius transformation, we can assume that the fixed points of $\beta$ are $0$ and $\infty$ and that the orbits of length $p$ are $$\{ u_{i,1}=\omega^i : i \in \mathbb{Z}_p\} \,\, \mbox{ and }\,\,\{ u_{i,k}=\lambda_k \omega^i :i \in \mathbb{Z}_p \} \, \mbox{ for } 2 \leqslant k \leqslant m,$$ where $\lambda_2, \ldots, \lambda_m$ are non-zero complex numbers as in the statement of the theorem (note that the conditions imposed on them guarantee that the orbits are disjoint).  Thus, after replacing $\omega$ by an appropriate power of it and after replacing $\lambda_k$ by $\omega^{u_k}\lambda_k$ for some $u_k$ if necessary,  we conclude that $S$ is isomorphic to the normalization of  \begin{equation} \label{model}y^q=\Pi_{i=0}^{p-1}(x-\omega^i)^{r^i}  \Pi_{k=2}^{m}\Pi_{i=0}^{p-1}(x- \lambda_k \omega^i)^{r^i}\end{equation}

\subsection*{Proof of statement (3)} Observe that with the previous  identification, the action of $P$ on $X$ is then given by $z \mapsto \beta(z)=\omega z$ and the regular covering map $X \to X/P$ is given by $z \mapsto z^p.$ Thus, after identifying $$S/G_{p,q} \cong X/P \cong \mathbb{P}^1,$$ the branch values of $S \to S/G_{p,q}$ are $\infty, 0$ marked with $p$ and $1, \lambda^p$ marked with $q.$ Note that the lift to $S$ of $\beta$ has the form $$(x,y) \mapsto (\omega x, \psi(x,y))$$where  $\psi(x,y)^q=f(\omega x)$ with $f(x)$ denoting the right side part of \eqref{model}. In fact $$\psi(x,y)=\varphi(x)y^r$$with $\varphi$ as in the theorem. Finally, it is a direct computation to verify that $$A(x,y)=(x, \xi y) \,\, \mbox{ and }\,\, B(x,y)=(\omega x , \varphi(x)y^r)$$satisfy $BAB^{-1}=A^r$ and therefore $\langle A, B \rangle \cong G_{p,q}$ as claimed.

\section{Proof of Theorem \ref{familia22}} \label{s4}
\subsection*{Proof of statement (1)}
Consider a Fuchsian group $\Delta$  of signature $s_{2,2}$  $$\Delta= \langle x_1,x_2,y_1,y_2 : x_1^p=x_2^p=y_1^q=y_2^q=x_1x_2y_1y_2=1\rangle$$and let $\theta: \Delta \to G_{p,q}$ be a ske representing an action of $G_{p,q}$ on $S \in \mathscr{C}_{2,2}$. Then $$\theta(x_1)=a^{l_1}b^{n_1} \,\, \mbox{ and } \theta(x_2)=a^{l_2}b^{n_2}$$for some $l_1, l_2 \in \mathbb{Z}_q$ and $n_1, n_2  \in \mathbb{Z}_p^*.$ After a suitable conjugation, we can suppose $l_2=0.$ Moreover, after applying an automorphism of $G$ of the form $a \mapsto a^i$ and  $b \mapsto b$ for some $i  \in \mathbb{Z}_q^*,$ we can assume that $\theta$ is given by $$\theta(x_1)=a^{l_1}b^{n_1}, \,\, \theta(x_2)=b^{n_2}, \,\, \theta(y_1)=a^{l_3} \,\, \mbox{ and }\,\, \theta(y_2)=a.$$The fact that $x_1x_2y_1y_2=1$ implies that $n_2=-n_1$ and $l_3=-l_1-1,$ and therefore $\theta$ is equivalent to the ske $\theta_{l,n}$ given by $$x_1 \mapsto a^{l}b^{n}, \,\, x_2 \mapsto b^{-n}, \,\, y_1 \mapsto a^{-l-1} \,\, \mbox{ and }\,\, y_2 \mapsto  a,$$for $l \in \mathbb{Z}_q$ and $n \in  \mathbb{Z}_p^*.$  Thus, there are at most $q(p-1)$ pairwise non-equivalent skes.

\subsection*{Proof of statement (2)}
As a consequence of \cite[Theorem 1]{singerman2}, the action of $G_{p,q}$ on $S$ might be only extended to an action of a group of order $2pq$ with signature $(0; 2,2,p,q)$. We claim that such extension is not possible. Indeed, if $G'$ is a group of order $2pq$ such that $G_{p,q} \leqslant G'$ then, by the  Schur-Zassenhaus theorem, we have that $$G' = G_{p,q} \rtimes C_2 \,\, \mbox{ where} \,\,\, C_2 = \langle t : t^2 = 1 \rangle.$$

As $\langle a \rangle$ contains all the elements of $G_{p,q}$ of order $q$ and as $t$ has order two, we observe that $tat^{-1}=a^{\epsilon}$ where $\epsilon = \pm 1$. We write $tbt^{-1}=a^{m}b^n$ for  $m \in \mathbb{Z}_q$ and $n \in \mathbb{Z}_p^*$ and observe that the fact that \begin{equation} \label{salida}b=t^2bt^{-2}=a^{m \epsilon}(a^mb^n)^n\end{equation}implies that $n=\pm1.$

\begin{enumerate}
\item Assume $\epsilon=1.$ The equality \eqref{salida} shows that  $m=0.$ Observe that  if $n=-1$  then $bab^{-1}=a^r$ implies $b^{-1}ab=a^{r};$ a contradiction. Consequently $n=1$ and $G'=G_{p,q} \times C_2.$ Note the $t$ is the unique involution of the group.
\item Assume $\epsilon = -1.$ If $n=-1$ then equality \eqref{salida} shows that $m=0,$ and $bab^{-1}=a^r$ implies $b^{-1}a^{-1}b=a^{-r};$ a contradiction.  It follows that $n=1.$  We write $b' := a^ub$ where $2u \equiv m \mbox{ mod }  q$ and notice that $tb't^{-1}=b'.$ Thus, we can assume $m=0$ and the involutions of the group are $ta^k$ where $k\in \mathbb{Z}_q.$\end{enumerate} 
The contradiction is obtained after noticing that if there were a ske from a Fuchsian group of signature $(0; 2,2,p,q)$ onto $G'$ then the group would contain an element of order $q$ and two involutions  whose product has order $p;$ however, in both cases, this is not possible.

\subsection*{Proof of statement (3)} Let $\Delta'$ be a  Fuchsian group of signature $(0; q, 2p, 2p)$  $$\Delta'=\langle z_1, z_2, z_3 : z_1^q=z_2^{2p}=z_3^{2p}=z_1z_2z_3=1\rangle$$and consider the group $G'$ of order $2pq$   $$G'=\langle a,c : a^q=c^{2p}=1, cac^{-1}=a^{-r}\rangle \cong C_q \rtimes C_{2p}.$$ Then the map $\Theta: \Delta' \to G'$ given by $$z_1 \mapsto a, \,\, z_2 \mapsto c \,\, \mbox{ and }\,\, z_3 \mapsto (ac)^{-1}$$is a ske and then the orbit space $Z=\mathbb{H}/\mbox{ker}(\Theta)$ is a compact Riemann surface  with a group of automorphisms isomorphic to $G'$ acting on it with signature $(0; q, 2p, 2p)$.

\s

We now proceed to prove that $Z \in \mathscr{C}_{2,2}.$ Define $$\hat{x}_1:=z_3^{-2} , \,\, \hat{x}_2:=z_1z_3^2z_1^{-1}, \,\, \hat{y}_1:=z_1 \,\, \mbox{ and } \,\, \hat{y}_2:=z_3^{-2}z_1^{-1}z_3^{2}$$and observe that $\hat{x}_1, \hat{x}_2$ have order $p,$  $\hat{y}_1, \hat{y}_2$ have order $q$ and the product of them equals 1. Note that, by letting $b=c^2,$ we have that  $$\Theta(\hat{x}_1)=a^{1-r}b, \,\, \Theta(\hat{x}_2)=a^{1+(r-2)r^{-1}}b^{-1}, \,\,\Theta(\hat{y}_1)=a \,\, \mbox{ and } \,\, \Theta(\hat{y}_2)=a^{-r}$$and therefore $$\Theta|_{\langle \hat{x}_1, \hat{x}_2, \hat{y}_1, \hat{y}_2 \rangle} : \langle \hat{x}_1, \hat{x}_2, \hat{y}_1, \hat{y}_2 \rangle \cong \Delta \to \langle a,b \rangle$$is equivalent to one of the skes $\theta_{l,n}$ as before. Consequently, $Z \in \mathscr{C}_{2,2}.$

\section{Proof of Theorem \ref{familia40}} \label{s6}
\subsection*{Proof of statement (1)} Consider a Fuchsian group $\Delta$  of signature $s_{0,4}$  $$\Delta= \langle x_1,x_2,x_4,x_4 : x_1^p=x_2^p=x_3^p=x_4^p=x_1x_2x_3x_4=1\rangle$$and let $\theta: \Delta \to  G_{p,q}$ be a ske representing an action of $G_{p,q}$ on $S \in \mathscr{C}_{4,0}$. Then $$\theta(x_i)=a^{l_i}b^{n_i} \,\, \mbox{ for } i=1,2,3,4$$where $l_i \in \mathbb{Z}_q$ and $n_i \in \mathbb{Z}_p^*.$ After a suitable conjugation, we can suppose $l_1=0.$ In addition, note that if each $l_i$ equals zero then $\theta$ is not surjective. Then,  after considering the braid transformations $\Phi_2$ or $\Phi_2 \circ \Phi_3$ if necessary, we can suppose $l_2 \neq 0.$ Now, we consider the automorphism of $G_{p,q}$ given by $a \mapsto a^{k_2}$ and  $b \mapsto b$ where $k_2l_2 \equiv 1 \mbox{ mod } q,$ to  assume $l_2=1.$ It follows that $\theta$ is given by $$\theta(x_1)=b^{n_1}, \,\, \theta(x_2)=ab^{n_2}, \,\, \theta(x_3)=a^{l_3}b^{n_3} \,\, \mbox{ and } \,\, \theta(x_4)=a^{l_4}b^{n_4}$$where, as $x_1x_2x_3x_4=1$, we have that $$n_1+n_2+n_3+n_4=0 \,\, \mbox{ and }\,\, 1+l_3r^{n_2}+l_4r^{n_2+n_3}=0.$$Thus, $\theta$ is equivalent to the ske $\theta_{n_1, n_2, n_3, l_3}$ given by $$x_1 \mapsto b^{n_1}, \,\, x_2 \mapsto ab^{n_2}, \,\, x_3 \mapsto a^{l_3}b^{n_3} \,\, \mbox{ and } \,\, x_4 \mapsto a^{-r^{-n_2-n_3}-l_3r^{-n_3}}b^{-n_1-n_2-n_3}$$ where $n_i \in \mathbb{Z}_p^*$ such that $n_1+n_2+n_3 \neq 0$ and $l_3 \in \mathbb{Z}_q.$  Then, there are at most $$q[(p-1)^3-(p-1)(p-2)]=q(p-1)(p^2-3p+3)$$pairwise non-equivalent skes.

\subsection*{Proof of statement (2)} By \cite[Theorem 1]{singerman2}, the action of $G_{p,q}$ on $S$ might be  extended to an action of a group $G'$ of order $2pq$ with signature $(0; 2,2,p,p)$ and this action, in turn, might be extended to only an action of a group of order $4pq$ with signature $(0; 2,2,2,p).$

\s

Note that $G'$ cannot be isomorphic to $G_{p,q} \times C_2$ since it has only one involution.
\s

{\bf Claim 1.} The action of $G_{p,q}$ on the surfaces lying in certain strata of $\mathscr{C}_{4,0}$ extends to an action of $$G' = \langle a,b,t : a^q=b^p=t^2=1, bab^{-1}=a^r, (ta)^2=[t,b]=1 \rangle \cong C_q \rtimes C_{2p}$$ with signature $(0; 2,2,p,p).$

\s

Let $\Delta'$ be a Fuchsian group of signature $(0; 2,2,p,p)$ $$\Delta' = \langle z_1, z_2, z_3, z_4 : z_1^2=z_2^2=z_3^p=z_4^p=z_1z_2z_3z_4=1\rangle,$$and let $\Theta: \Delta' \to G'$ be the ske defined by $$\Theta(z_1)=ta^{L}, \,\, \Theta(z_2)=ta^{L-1}, \,\,  \Theta(z_3)=ab^{N} \,\, \mbox{ and }\,\, \Theta(z_4)=b^{-N},$$where $L \in \mathbb{Z}_q$ and $N \in \mathbb{Z}_p^*.$ Set$$\hat{x}_1:=z_3, \,\, 
\hat{x}_2:= z_4, \,\,
\hat{x}_3:=z_1z_3z_1 \,\, \mbox{ and } \,\, 
\hat{x}_4:=z_1z_4z_1$$and note that they generate a Fuchsian group isomorphic to $\Delta.$ The claim follows after noticing that the restriction$$\Theta|_{\langle \hat{x}_1, \hat{x}_2, \hat{x}_3, \hat{x}_4 \rangle} : \Delta \to \langle a,b \rangle \cong G_{p,q}$$is equivalent to some ske $\theta_{n_1, n_2, n_3, l_3}$ as before, with $n_2=-n_1=N$.

\s

{\bf Claim 2.} The action of $G_{p,q}$ on each $S$ does not extend to any action of a group of order $4pq$ with signature $(0; 2,2,2,p).$

\s

Assume that the action of $G_{p,q}$ on $S$ extends to an action of a group $G''$ of order $4pq$ with signature $(0; 2,2,2,p).$ Then $G''$ is isomorphic to either $$ (C_q \rtimes C_p) \rtimes C_4  \,\, \mbox{ or }\,\, (C_q \rtimes C_p) \rtimes C_2^2 .$$

\s

In the former case, if $C_4 = \langle t : t^4 =1 \rangle$ then $$tat^{-1}=a^{\epsilon} \,\, \mbox{ and } \,\, tbt^{-1}=a^nb$$where $n \in \mathbb{Z}_q$ and $\epsilon$ equals 1, $-1$ or a primitive fourth root of unity in the field of $q$ elements. If $\epsilon =1$ then the product is direct, if $\epsilon = -1$ then $t^2$ is the unique involution and if $\epsilon$ is a primitive fourth root of unity in the field of $q$ elements then the involutions are $t^2a^k$ where $k \in \mathbb{Z}_q.$ 

\s

In the latter case, if $C_2^2 = \langle t, u : t^2=u^2=(tu)^2 =1 \rangle$ then $$tat=a^{\epsilon_1}, \,\, tbt=a^{n_1}b, \,\,uau=a^{\epsilon_2} \,\, \mbox{ and }\,\,  ubu=a^{n_2}b$$where $n_i \in \mathbb{Z}_q$ and $\epsilon_i =\pm 1.$  If $\epsilon_1=\epsilon_2=1$ then the product is direct, if $\epsilon_1=\epsilon_2=-1$ then the involutions are $ta^k$ and $ua^k$ where $k \in \mathbb{Z}_q,$ and if $\epsilon_1=-\epsilon_2=1$ or $\epsilon_1=-\epsilon_2=-1$ then the involutions are $ta^k$ or $ua^k$ respectively, where $k \in \mathbb{Z}_q.$ 

\s

The proof of the claim follows after noticing that, in each case, the group $G''$ cannot be generated by three involutions.

\subsection*{Proof of statement (3)} Let $\Delta''$ be a  Fuchsian group of signature $(0; 2, 2p, 2p)$ $$\Delta''=\langle z_1, z_2, z_3 : z_1^2=z_2^{2p}=z_3^{2p}=z_1z_2z_3=1\rangle$$and consider the group $G''=(C_q \rtimes C_{2p}) \rtimes C_2$ presented as $$\langle a, c, t: a^q=c^{2p}=t^2=1, cac^{-1}=a^{-r}, (ta)^2 =[t,c]=1\rangle.$$ The map $\Theta: \Delta'' \to G''$ given by $$z_1 \mapsto ta, \,\, z_2 \mapsto a^{-1}c \,\, \mbox{ and }\,\, z_3 \mapsto c^{-1}t$$is a ske and therefore $V=\mathbb{H}/\mbox{ker}(\Theta)$ is a compact Riemann surface with a group of automorphisms isomorphic to $G''$ acting on it with signature $(0; 2, 2p, 2p)$. 

\s

 If we define $$\hat{x}_1:=z_3^{2} , \,\, \hat{x}_2:=z_1z_2^2z_1, \,\, \hat{x}_3:=z_1z_3^2z_1 \,\, \mbox{ and } \,\, \hat{x}_4:=z_2^2$$then they generate a Fuchsian group isomorphic to $\Delta$ and the restriction $$ \Theta|_{\langle \hat{x}_1, \hat{x}_2, \hat{x}_3, \hat{x}_4 \rangle} : \Delta \to \langle a,b:=c^2 \rangle \cong G_{p,q}$$is equivalent to some ske $\theta_{n_1, n_2, n_3, l_3}$ as before. Hence, $V \in \mathscr{C}_{4,0}.$

\subsection*{Proof of statement (4)} Consider the subgroups $N=\langle a \rangle$ and $K=\langle a, t \rangle$ of $G'$ and the quotient $J=K/N \cong C_2$. Note that $K$  is isomorphic to the dihedral group of order $2q$. The regular covering map $S \to S/K$ ramifies over exactly $2p$ values  marked with two and, consequently, the quotient $S/K$ has genus zero. It is clear that $N$ acts on $S$ without fixed points and therefore the genus of the quotient $X=S/N$ equals $\gamma=p-1.$ Note that $X$ admits the action of $J$ and the corresponding two-fold regular covering map $$X \to  X/J \cong S/K \cong \mathbb{P}^1$$ramifies over $2p=2\gamma+2$ values, showing that $X$ is hyperelliptic. We observe that, as the group $P=G'/K \cong C_p$ acts on $X/J$ with signature $(0; p,p),$ we can suppose that $P \cong \langle \beta \rangle$ where $\beta(z):=\omega z,$ and that the branch values of $X \to X/J$ are  $$\{\omega^k : k \in \mathbb{Z}_p\} \,\, \mbox{ and } \{\mu \omega^k : k \in \mathbb{Z}_p\} $$where $\omega$ is a primitive $p$-th root of unity and $\mu$ is a non-zero complex number such that $\mu^p \neq 1.$  The result follows by arguing as in the proof of Theorem \ref{modelo}.
\section{Proof of Theorem \ref{familia31}} \label{s5}

The proof 
 is similar to the ones of Theorems \ref{familia22} and \ref{familia40}; so, we avoid some details.
 
\subsection*{Proof of statement (1)} Consider a Fuchsian group $\Delta$ of signature $s_{3,1}$  $$\Delta= \langle x_1,x_2,x_3,y_1 : x_1^p=x_2^p=x_3^p=y_1^q=x_1x_2x_3y_1=1\rangle$$and let $\theta: \Delta \to G_{p,q}$ be a ske representing an action of $G_{p,q}$ on $S \in \mathscr{C}_{3,1}$. Up to equivalence, we can suppose  that $\theta$ agrees with the ske $\theta_{l,n_1, n_2}$ given by $$x_1 \mapsto a^{-1-lr^{n_1}}b^{n_1}, \,\, x_2 \mapsto a^{l}b^{n_2}, \,\, x_3 \mapsto b^{-n_1-n_2} \,\, \mbox{ and }\,\, y_1 \mapsto a,$$for some  $l \in \mathbb{Z}_q$ and $n_i \in \mathbb{Z}_p^*$ such that $n_1 \neq -n_2.$  Thus, there are at most $q(p-1)(p-2)$ pairwise non-equivalent skes.

\subsection*{Proof of statement (2)}
It follows from the fact that, by \cite[Theorem 1]{singerman2}, the  signature $(0; p,p,p,q)$ is maximal.

\subsection*{Proof of statement (3)} Let $\Delta'$ be a  Fuchsian group of signature $(0; p, 2p, 2q)$ $$\Delta'=\langle z_1, z_2, z_3 : z_1^p=z_2^{2p}=z_3^{2q}=z_1z_2z_3=1\rangle$$and consider the direct product $G'=G_{p,q} \times C_2$ with $C_2=\langle t : t^2 = 1 \rangle.$ The map $$\Theta: \Delta' \to G' \,\, \mbox{ given by } z_1 \mapsto a^{-1}b, \,\, z_2 \mapsto b^{-1}t \,\, \mbox{ and }\,\, z_3 \mapsto ta$$is a ske and then $W=\mathbb{H}/\mbox{ker}(\Theta)$ is a  Riemann surface with a group of automorphisms isomorphic to $G'$ acting on it with signature $(0; p, 2p, 2q)$. To see that $W \in \mathscr{C}_{3,1},$ define $$\hat{x}_1:=z_1 , \,\, \hat{x}_2:=z_2^2, \,\, \hat{x}_3:=z_3z_1z_3^{-1} \,\, \mbox{ and } \,\, \hat{y}_1:=z_3^2$$and notice that they generate a Fuchsian group isomorphic to $\Delta$ and that the restriction of $\Theta$ to it is equivalent to some ske $\theta_{l,n_1, n_2}$ as before.

\section{Proof of Theorems \ref{jaco} and \ref{prop3}} \label{s7}
\subsection*{Rational irreducible representations} We refer to \cite{Serre} for basic background concerning representations of groups.
Set $ \omega_t:=\mbox{exp}(\tfrac{2 \pi i}{t})$ where $i^2=-1.$ The group $G_{p,q}$ has, up to equivalence, $p$ complex irreducible representations of  degree 1, given by $$\chi_l : a \mapsto 1 \,\, \mbox{ and } \,\, b \mapsto\omega_{p}^l \,\, \mbox{ for }l \in \mathbb{Z}_p.$$They give rise to two non-equivalent rational irreducible representations of the group; namely, the trivial one $W_0=\chi_0$ and the direct sum of the non-trivial ones $$W_1=\chi_1 \oplus \cdots \oplus \chi_{p-1}.$$

Consider the equivalence relation $R$ on $\mathbb{Z}_q^*$ given by $$uRv \iff u=r^nv \, \mbox{ for some } n \in \mathbb{Z}_p$$where $r$ is primitive $p$-th root of unity  in the field of $q$ elements. Let $d=\tfrac{q-1}{p}$ and let $\{k_1, \ldots, k_{d}\}$ be a maximal collection of representatives of $R$. Then $G_{p,q}$ has, up to equivalence,
$d$ complex irreducible representations of degree $p$$$\psi_j: a \mapsto \mbox{diag}(\omega_q^{k_j}, \omega_q^{k_jr}, \omega_q^{k_jr^2}, \ldots, \omega_q^{k_jr^{p-1}}) \,\, \mbox{ and } \,\, b \mapsto \left( \begin{smallmatrix}
0 & 1 & 0 & \cdots & 0 \\
0 & 0 & 1 & \cdots & 0 \\
\, & \, & \, & \ddots & \, \\
0 & 0 & 0 & \cdots & 1 \\
1 & 0 & 0 & \cdots & 0 \\
\end{smallmatrix} \right)
$$for $ j \in \{1, \ldots, d\}.$ The direct sum of these representations yields a rational irreducible representation of the group; namely $$W_2= \psi_1 \oplus \cdots \oplus \psi_d.$$

Observe that the set $\{W_0, W_1, W_2\}$ is a maximal set of pairwise non-equivalent rational irreducible representations of $G_{p,q}.$ In addition, $\chi_1$ and $\psi_1$ have Schur index one  and their character fields have degree $p-1$ and $(q-1)/p$ over $\mathbb{Q}$ respectively.

\s

\subsection*{Proof of Theorem \ref{jaco}}Let $S$ be a compact Riemann surface endowed with a group of automorphisms isomorphic to $G_{p,q}.$ Then, as explained in Subsection \S\ref{ii3},   the information concerning the rational irreducible representations  described above, allows us to obtain that the group algebra decomposition of $JS$ with respect to $G_{p,q}$ is \begin{equation}\label{gad}JS \sim B_0 \times B_1 \times B_2^p\end{equation}where the factor $B_l$ is associated to the representation $W_l$ (and, in turn, $W_1, W_2$ and $W_3$ are associated to $\chi_0, \chi_1$ and $\psi_1$ respectively). Note that $B_0$ is isogenous to the Jacobian variety of $S/G_{p,q}.$
\s

We now assume that $S \in \mathscr{C}_{n,m}$ and that  the action of $G_{p,q}$ on $S$ is represented by the ske $\theta: \Delta \to G_{p,q},$ with $\Delta$ presented as in the proof of Theorem \ref{tepi}. Observe that, independently of $j$ and independently of the choice of the ske $\theta$, we have that $$\langle \theta(y_j) \rangle = \langle a \rangle \,\, \mbox{ and } \,\, \langle \theta(x_j) \rangle \sim_c \langle b \rangle$$where $\sim_c$ stands for conjugation. It follows that the dimension of the fixed  subspaces of $\chi_1$ and $\psi_1$ under the action of the corresponding isotropy groups $\langle \theta(x_j) \rangle$ and $\langle \theta(y_j) \rangle$ do not depend on $\theta,$ and are given by $$\mbox{dim} (\chi_1^{\langle \theta(x_j) \rangle})=0 \, \mbox{ and }\, \mbox{dim} (\chi_1^{\langle \theta(y_j) \rangle})=1$$and $$\mbox{dim} (\psi_1^{\langle \theta(x_j) \rangle})=1 \, \mbox{ and }\, \mbox{dim} (\psi_1^{\langle \theta(y_j) \rangle})=0$$for each $j.$ We now apply the equation \eqref{ooii} to conclude that the dimension of the factors $B_1$ and $B_2$ in \eqref{gad} are given by $$\dim(B_1)=(p-1)[-1+\tfrac{1}{2}(n(1-0)+m(1-1))]=(n-2)(\tfrac{p-1}{2})$$ $$\dim(B_2)=\tfrac{q-1}{p}[-p+\tfrac{1}{2}(n(p-1)+m(p-0))]=\tfrac{q-1}{p}[p(\tfrac{m}{2}-1)+n(\tfrac{p-1}{2})],$$as claimed. Clearly $B_0=0$ since $S/G_{p,q} \cong \mathbb{P}^1.$

\subsection*{Proof of Theorem \ref{prop3}} Let $S$ be a compact Riemann surface lying in the family $\mathscr{C}_{n,m}$ and consider the group algebra decomposition of $JS$ with respect to $G_{p,q}$ \begin{equation*} 
JS \sim B_1 \times B_2^p
\end{equation*}obtained in Theorem \ref{jaco}. Following \cite[Proposition 5.2]{cr}, if  $K$ is a subgroup of $G_{p,q}$ then   \begin{equation} \label{my}J(S/K) \sim  B_{1}^{{n}_1^K}  \times B_{2}^{n_2^K}
\end{equation}where $n_1^K$ and $n_2^K$ are the dimension of the vector subspaces of $\chi_1$ and $\psi_1$ fixed under the action of $K.$ Now, if we consider \eqref{my} with $K=N$ and $K=H$ we obtain $$JX=J(S/N)\sim B_1^1 \times B_2^0 \,\, \mbox{ and }\,\, JY=J(S/H) \sim B_1^0 \times B_2^1$$respectively, and the result follows.

\s

\subsection*{Acknowledgements} The author is    grateful to the referee for valuable comments and suggestions.

\end{document}